\documentclass{gtpart}     

\usepackage{etex}
\usepackage{indentfirst}
\usepackage[below]{placeins} 

\usepackage{xcolor} 

\usepackage{amsmath,amssymb,amsfonts,amsthm,latexsym}
\usepackage{amsfonts}
\usepackage[utf8]{inputenc} \usepackage{ae,aecompl}     
\usepackage{amscd}
\usepackage{mathrsfs}       
\usepackage{fancyhdr}
\usepackage{graphicx,psfrag}
\usepackage{epsf} 
\usepackage[all]{xy}
\newtheorem{thm}{Theorem}[section]
\newtheorem{innercustomthm}{Theorem}
\newenvironment{customthm}[1]
  {\renewcommand\theinnercustomthm{#1}\innercustomthm}
  {\endinnercustomthm}

\newtheorem{lem}[thm]{Lemma}

\newtheorem{cor}[thm]{Corollary}
\newtheorem{pro}[thm]{Proposition}
\newtheorem{defn}[thm]{Definition}

\theoremstyle{definition}

\newcommand {\Ta}{\mathbb{ T}_{\alpha}}
\newcommand {\Tb}{\mathbb{ T}_{\beta}}

\newcommand{\hfhat}{\widehat{HF}}

\newcommand{\cfkhat}{\widehat{CFK}}
\newcommand{\cfkinfty}{CFK^{\infty}}

\newcommand{\hfk}{\widehat{HFK}}
\newcommand{\hfkhat}{\widehat{HFK}}

\newcommand{\x}{{\bf x}}

\newcommand{\s}{\mathfrak{s}}

\newcommand{\rank}{\text{rank}}

\newcommand{\spinc}{\mathrm{ Spin}^c }

\title{Knot Complement Problem for L-space $\Z HS^3$}
\author{Huygens C. Ravelomanana}
\givenname{Huygens C.}
\surname{Ravelomanana}
\address{Department of Mathematics\\
452 Boyd Graduate Studies\\
University of Georgia\\
Athens, GA 30602, US}
\email{huygens@cirget.ca}
\urladdr{}

\arxivreference{arXiv:1505.00239v3}  
\arxivpassword{}   
\volumenumber{}
\issuenumber{}
\publicationyear{}
\papernumber{}
\startpage{}
\endpage{}
\doi{}
\MR{}
\Zbl{}
\received{}
\revised{}
\accepted{}
\published{}
\publishedonline{}
\proposed{}
\seconded{}
\corresponding{}
\editor{}
\version{}

\begin{document}

\begin{abstract}
In this paper we look at the knot complement problem for L-space $\Z$-homology spheres. We show that an L-space $\Z$-homology sphere $Y$ cannot be obtained as a non-trivial surgery along a knot $K\subset Y$. As a consequence, we prove that  knots in an L-space $\Z$-homology sphere are determined by their complements. 
\end{abstract}
\maketitle

\section{Introduction}

 A knot $K$ in a 3-manifold $Y$ is \textit{determined by its complement} if the existence of a homeomorphism between $Y\setminus K$ and $Y\setminus K'$ for some other knot $K$, implies the existence of homeomorphism between the pair $(Y,K)$ and $(Y,K')$. Here we do not expect these homeomorphisms to be orientation preserving. The \textit{``knot complement problem"} is the problem of knowing if a given knot is determined by its complement. Note that in this context two knots are equivalent if there is an homeomorphism of the ambient manifold which send one to the other and the homeomorphism need not be orientation preserving. Gordon and Luecke \cite{Gordon-Luecke} have proved that non-trivial knots in $S^3$ and $S^2\times S^1$ are determined by their complements. D. Matignon has proved in \cite{Matignon-knot} that, if one only considers orientation preserving homeomorphisms, then  non-trivial non-hyperbolic knots are determined by their complements in closed, atoroidal and irreducible Seifert fibred 3-manifolds; except the axes in $L(p,q)$ when $q^2\equiv \pm 1\ [\text{mod}\ p]$. 
In general there are some knots which are not determined by their complements. In \cite{Rong} Y. Rong classified all such Seifert fibered knots in closed 3-manifolds other than lens spaces. In \cite{Matignon-hyp-knot} D. Matignon give examples of hyperbolic knots in lens spaces which are not determined by their complements. Here we are interested in a problem related to the knot complement problem for L-space $\Z$-homology spheres. A example of this are knots in the Poincar\'e sphere $\Sigma(2,3,5)$. More precisely, we prove the following theorems and corollaries.

\begin{thm}\label{my-theorem-1-general}
Let $K$ be a knot in an oriented L-space $\Z$-homology sphere $Y$ and let $r\in \Q$. The result of an $r$-surgery along $K$ is never  homeomorphic to $Y$.
\end{thm}

From this theorem we can answer the  oriented knot complement problem for L-space $\Z$-homology spheres.

\begin{thm}\label{my-theorem-3}
Knots in L-space $\Z$-homology spheres are determined by their oriented complements.
\end{thm}

In particular this is true for the Poincar\'e sphere $\Sigma(2,3,5)$ which is the only known irreducible L-space $\Z$-homology sphere apart from $S^3$. Note also that $\Sigma(2,3,5)$ does not admit orientation reversing homeomorphisms.

\begin{cor}\label{my-theorem-1}
Let $K$ be a knot in $\Sigma(2,3,5)$ and let $r\in \Q$. The result of an $r$-surgery along $K$ is never  homeomorphic to $\Sigma(2,3,5)$.
\end{cor}

\begin{cor}\label{my-theorem-2}
Knots in $\Sigma(2,3,5)$ are determined by their complements.
\end{cor}


In the earlier version of the paper Theorem \ref{my-theorem-1-general} and Theorem \ref{my-theorem-3}  was stated only for orientation preserving homeomorphisms. Fyodor Gainullin \cite{Gainullin}  independently proved that more generally null-homologous knots in a rational homology L-space are determined by their oriented complements. Note that Lens spaces are rational homology L-spaces but the axes in Lens spaces are not null-homologous and they are knot determined by their complements by \cite{Matignon-knot}. We prove here the special case of L-space $\Z$-homology spheres.  The main part of the proof here, Lemma \ref{my-lemma-sigma235-1}, was a lemma in the author Ph.D. thesis and is a direct consequence of a previous work of J. Rasmussen. The main theorem is obtained by using properties of genus 1 fibred knots together with Lemma \ref{my-lemma-sigma235-1}. 

\subsection*{Organization}
The paper is organized as follows. In section~\ref{section: preliminaries} we give some preliminaries. In section~\ref{section: proof of main theorems} we give the proof of the main theorems.

\subsection*{Acknowledgment.}
I would like to thank my supervisor Steven Boyer for his support. I would like to thank  Andr\'as Stipsicz for some helpful discussions. I especially thank Kenneth Baker for pointing out the uniqueness of genus one fibred knot in the Poincar\'e sphere and for helpful conversations. I would like to thank Daryl Cooper for comments about the orientation. This paper resulted from works on the more general cosmetic surgery problem carried out while the author was a graduate student at UQ\`AM and CIRGET in Montr\'eal. This version was updated when the author was at the University of Georgia.

\section{Preliminaries}\label{section: preliminaries}

\subsubsection*{Heegaard Floer homologies.}

Heegaard Floer homology is an invariant for closed oriented three manifolds $Y$.  The invariant, denoted $\hfhat$ is the homology of a chain complex defined from an Heegaard splitting of $Y$ but also admits some combinatorial definition. Ozsv\'ath-Szab\'o in \cite{Holo-disks-knot} and Rasmussen in \cite{Rasmussen-Knot} defined a related invariant for null-homologous knots $K$ in $Y$, taking the form of an induced filtration on the Heegaard Floer complex of $Y$.  The filtered chain homotopy type of this complex is a knot invariant, the ``knot Floer homology''. 

Knot Floer homology then associates to a null-homologous knot $K$ in $Y$ a 
$\Z\oplus\Z$--filtered complex $\cfkinfty(Y,K)$, generated over $\Z$ by elements $[\x;i,j]$ where $x\in \Ta\cap \Tb$  and $(i,j)\in \Z\oplus\Z$.  The set $\Ta\cap \Tb$ is the intersection of two half dimensional totally real submanifolds in a symplectic manifold obtained from an Heegaard splitting of $Y$.  Let $\mathfrak{G}$  be the set of generator of $\cfkinfty(Y,K)$, the $\Z\oplus\Z$--filtration on the complex is  given by the map $\mathcal{F}:\mathfrak{G} \longrightarrow\Z\oplus \Z$ defined by $\mathcal{F}([\x;i,j])=(i,j)$. The complex $\cfkhat(Y,K)$ is a subcomplex of a quotient complex of $\cfkinfty(Y,K)$ and is equipped with a $\Z$-filtration. The homology of  $\cfkhat(Y,K)$ is denoted  $\hfkhat(Y,K)$ and its Euler characteristic is the normalized Alexander polynomial $\Delta_K(T)$ of the knot $K$ with the $\Z$-filtration corresponding to the exponent of $T$. From more details on the subject we refer to [\cite{Holo-disks-closed-3mfd}, \cite{Applications}, \cite{Holo-disks-knot}, \cite{Intro-HFH}, \cite{Lectures-HFH}, \cite{Rasmussen-Knot}].

\subsubsection*{$L$-space.}
The Heegaard Floer homology $\hfhat(Y)$ has the property that 
$$\rank \hfhat(Y) \geq |H_1(Y;\Z)|.$$
The manifolds which satisfy the equality form an  important  class of manifold in 3-manifold theory.
\begin{defn} \label{definition:L-space}
A rational homology sphere $Y$ is called an L-space if $$\rank \hfhat(Y) = |H_1(Y;\Z)|.$$
\end{defn}
Being an L-space is also equivalent to $\hfhat(Y,\s) \cong \Z$ for each $\spinc$ structure $\s$ on $Y$.

L-spaces are Heegaard Floer analogues of lens-spaces. In particular every Lens space $L(p,q)$ with $p\neq 0$ is an L-space. They also include double branched cover of non-split alternating links.

A knot $K \subset Y$ is said to admit an $L$-space surgery if for some rational number $r$, $Y_K(r)$ is an $L$-space. We will us the fact that such knot has a very special knot Floer homology.

\subsubsection*{Irreducibility of $Y\setminus K$}
Before continuing further let us discuss the  irreducibility of $Y\setminus K$. The knot complement $Y\setminus K$  is irreducible if and only if it does not lie in a ball. Indeed if $Y\setminus K$ was reducible then there are two 3-manifolds $M$ and $N$   such that $M$ is distinct from $S^3$, $K\subset N$ and $Y= M \sharp  \left(N\setminus K \right)$. However meridional surgery on $K$  yields $Y$
 again so $Y_K(\infty)=Y$. Thus
$$Y_K(\infty)=Y=M \sharp N_K(\infty).$$
But $Y$ is irreducible and $M$ is distinct from $S^3$. Therefore $N_K(\infty)\cong S^3$ and it follows that $N\cong S^3$. Therefore $K$ lies in a ball. Conversely if a non-trivial knot $K$ in $Y$ lies in ball, then its the complement is obviously reducible.

The problem for a non-trivial knot which lies in a ball is then equivalent to the  problem for knots in $S^3$. Since it is known \cite{Gordon-Luecke} that no non-trivial surgery on $S^3$ can reproduce $S^3$, we can assume that our knot $K$ does not lie in a ball. That is we can assume $Y\setminus K$ is irreducible. 

\section{Proof of main theorems}\label{section: proof of main theorems}
We need the following  characterization of $\hfk(Y,K)$ for a $K$ is a knot in an integral homology $L$-space $Y$. It was proved in [\cite{OSzLensSpace} theorem 1.2] for the case of $S^3$ and stated in [\cite{Wu} proposition 3.7] for the more general case of $L$-space homology sphere. 

\begin{pro}\label{theorem-HFK-for-L-space-surgery}
Let $K$ be a knot in a $L$-space $\Z$-homology sphere $Y$.  If $Y_K(r)$ is an $L$-space for some rational number $r$, then there is an increasing sequence of integers $n_{-k}<\cdots<n_k$
such that $n_i=-n_{-i}$,  and $\hfk(K,j)=0$ unless $j=n_i$ for some $i$, in which case
$\hfk(K,j)\cong \Z$.
\end{pro}

An immediate corollary [\cite{Wu} corollary 3.8] is a simplified expression for the Alexander polynomials of such knots.

\begin{cor}\label{corollary-Alexander_poly}
Let $K$ be a knot in a $L$-space $\Z$-homology sphere. If $K$ admits an $L$-space surgery, then the Alexander polynomial of $K$ has the form
$$\Delta_K(T) = (-1)^k+ \sum_{j=1}^k(-1)^{k-j} (T^{n_j}+T^{-n_j}),$$
for some increasing sequence of positive integers $0<n_1<n_2<\cdots<n_k$.
\end{cor}

Our key lemma is the following which is about a characterisation of the Alexander polynomial of knots
in L-space $\Z$-homology spheres which admit L-space $\Z$-homology sphere surgery.

\begin{lem}\label{my-lemma-sigma235-1}
 Let $Z$ be an L-space $\Z$-homology sphere and let  $K$ be any non-trivial knot in $Z$. If $K$ admits an L-space $\Z$-homology sphere non-trivial surgery, then $\Delta_K(T)=  T^{-1} - 1 + T$, $K$ has genus $1$ and  $\hfkhat(K)\cong \Z^3$.
\end{lem}
Here the Alexander polynomial is normalized so that $\Delta_K(1)=1$ and $\Delta_K(T)=\Delta_K(T^{-1})$.
\begin{proof}
 Rasmussen proved in [\cite{Lens-Rasmussen} Proposition 4.5] that:

If $Z$ is an $L$-space with $H_1(Z)=\Z/p\Z$, and $K\subset Z$ is a \textit{primitive knot} (\textit{i.e} $K$ generates $H_1(Z)$) with  a homology sphere non-trivial surgery $X$. Then $X$ is an L-space if and only if one of the following condition  holds:
\begin{enumerate}
\item $\hfk (K) \cong \Z^p$ and width $\hfk (K) < 2p$.
\item $\hfk (K) \cong \Z^{p+2}$,  width $\hfk (K) = 2p$.
\end{enumerate}
Here width $\hfk (K) $  is the difference $Max - Min$, where $Max$ is the maximum value of $j$ for which $\hfk(K,j)$ is nontrivial and  $Min$ is the minimum value.

On the other hand we also have a formula for the genus of $K$, $g(K)= (\text{width}\ \hfk (K) -p+1)/2 $. See \cite{Lens-Rasmussen} theorem 4.3. 

In our case $p=1$ so all the hypothesis are satisfied and we have either 
\begin{enumerate}
\item $\hfk (K) \cong \Z$ and width $\hfk (K) < 2$.
\item $\hfk (K) \cong \Z^3$,  width $\hfk (K) = 2$.
\end{enumerate}

The first case implies that width $\hfk (K) =0$ then $g(K)=0$. Therefore we are left with the second case, we can then compute the Euler characteristic of $\hfk (K)$ to obtain the symmetrized Alexander polynomial of $K$ using the formula 
$$\Delta_K(T)= \sum_{i,j}\  (-1)^i \ \dim \hfk_i (K,j)\ T^j,$$
we obtain
$$\Delta_K(T)= a_{j_0}\ T^{j_0}+a_{j_0+1}\ T^{j_0+1}+a_{j_0+2}\ T^{j_0+2},$$
for some $j_0$, using the fact that $\Delta_K(T)=\Delta_K(T^{-1})$ we get  $j_0=-1$. Since  width $\hfk (K) = 2$ we have $a_{-1}=a_{1}=\pm 1$. so
$$\Delta_K(T)= \pm  T^{-1} + a_0  \pm T.$$
On the other hand From corollary \ref{corollary-Alexander_poly} $$\Delta_K(T) = (-1)^k+ \sum_{j=1}^k(-1)^{k-j} (T^{n_j}+T^{-n_j}),$$
for some increasing sequence of positive integers $0<n_1<n_2<...<n_k$. Therefore $j=1, k=1,\ n_j=1$ and
$$\Delta_K(T)=  T^{-1} - 1 + T.$$
Now computing the second derivative gives $\Delta_K''(1)=  2$. Finally since width $\hfk (K) = 2$ and $g(K)=\max\{k\ |\  \hfk_* (K,k) \neq 0\},$ we must have $g(K)=1$.
\end{proof}

\begin{customthm}{\ref{my-theorem-1-general}}
Let $K$ be a knot in an oriented L-space $\Z$-homology sphere $Y$ and let $r\in \Q$. The result of an $r$-surgery along $K$ is never homeomorphic to $Y$.
\end{customthm}

\begin{proof}
Since  $Y$ is an L-space $\Z$-homology sphere, lemma \ref{my-lemma-sigma235-1} implies that $K$ is a genus one knot, $\Delta_K(T)=  T^{-1} - 1 + T$ and $\hfkhat(K)\cong \Z^3$. Since we can assume $Y\setminus K$ is irreducible, the result of Yi Ni about fibred knot [\cite{Ni-fibred-knot} Theorem 1.1] applies. Therefore $K$ is a genus one fibred knot since $\hfkhat(K)\cong \Z^3$.  By chapter 5 of \cite{Burde-Zieschang}, since $\Delta_K(T)=  T^{-1} - 1 + T$, the monodromy of the fibered knot $K$ is the monodromy of the trefoil knot so $Y_K$ is the trefoil exterior. It then follows that $Y_K(r)$ is cannot be homeomorphic to $Y$.

\end{proof}

From this theorem we can answer the  oriented knot complement problem for knots in  L-space $\Z$-homology spheres.  

\begin{customthm}{\ref{my-theorem-3}}
Knots in L-space $\Z$-homology spheres are determined by their oriented complements.
\end{customthm}
\begin{proof}
Let $K$ and $K'$ be two non-trivial knots in an L-space $\Z$-homology sphere $Y$, let us denote $V$ and $V'$ their complements with the induced orientations. Suppose there is an orientation preserving homeomorphism $f:V\to V'$. Let $\mu_K$, respectively $\mu_{K'}$, be the meridional slope of $K$, respectively $K'$, and let $r=f(\mu_K)$. The oriented manifold $V'(r)$ is orientation preserving homeomorphic to $Y$ and therefore by Theorem \ref{my-theorem-1-general} $r=\pm \mu_{K'}$. It follows that we can extend $f$ to an orientation preserving homeomorphism between $\left(Y, K\right)$ and $\left(Y, K'\right)$.
\end{proof}

For the special case of $\Sigma(2,3,5)$ Kenneth Baker (\cite{Ken} Lemma 13) showed that there is a unique genus one fibred knot in $\Sigma(2,3,5)$ which is the surgery dual to the negative trefoil. This gives another proof that knots in  $\Sigma(2,3,5)$ are determined by their complements.

\bibliographystyle{gtart}
\bibliography{biblio-knot-complement-for-ZHS}

\begin{thebibliography}{}
\providecommand\bibmarginpar{\leavevmode\marginpar}
\def\urlstyle#1{{\tt #1}}

\bibitem{Ken}
\textbf{K\,L Baker}, \emph{The {P}oincar\'e homology sphere, lens space
  surgeries, and some knots with tunnel number two} (2015)
  \xox{arXiv}{1504.06682v1}

\bibitem{Burde-Zieschang}
\textbf{G Burder}, \textbf{H Zieschang}, \emph{Knots}, De Gruyter studies in
  mathematics, W. de Gruyter, Berlin (1985)

\bibitem{Gainullin}
\textbf{F Gainullin}, \emph{Heegaard Floer homology and knots determined by
  their complements} (2015) \xox{arXiv}{1504.06180}

\bibitem{Gordon-Luecke}
\textbf{C\,M Gordon}, \textbf{J Luecke}, \emph{Knots are Determined by Their
  Complements}, Journal of the American Mathematical Society 2 (1989) 371--415

\bibitem{Matignon-hyp-knot}
\textbf{D Matignon}, \emph{Problem of hyperbolic knot complement in lens
  spaces,}Preprint

\bibitem{Matignon-knot}
\textbf{D Matignon}, \emph{On the knot complement problem for non-hyperbolic
  knots}, Topology Appl. 157 (2010) 1900--1925

\bibitem{Ni-fibred-knot}
\textbf{Y Ni}, \emph{Knot {F}loer homology detects fibred knots}, Invent. Math.
  170 (2007) 577--608

\bibitem{Holo-disks-knot}
\textbf{P Ozsv{\'a}th}, \textbf{Z Szab{\'o}}, \emph{Holomorphic disks and knot
  invariants}, Adv. Math. 186 (2004) 58--116

\bibitem{Applications}
\textbf{P Ozsv{\'a}th}, \textbf{Z Szab{\'o}}, \emph{Holomorphic Disks and
  Three-Manifold Invariants: Properties and Applications}, Annals of
  Mathematics 159 (2004) pp. 1159--1245

\bibitem{Holo-disks-closed-3mfd}
\textbf{P Ozsv{\'a}th}, \textbf{Z Szab{\'o}}, \emph{Holomorphic disks and
  topological invariants for closed three-manifolds}, Ann. of Math. (2) 159
  (2004) 1027--1158

\bibitem{OSzLensSpace}
\textbf{P Ozsv{\'a}th}, \textbf{Z Szab{\'o}}, \emph{On knot {F}loer homology
  and lens space surgeries}, Topology 44 (2005) 1281--1300

\bibitem{Intro-HFH}
\textbf{P Ozsv{\'a}th}, \textbf{Z Szab{\'o}}, \emph{Introduction to {H}eegaard
  {F}loer theory}, from: ``Floer Homology, Gauge Theory, and Low Dimensional
  Topology'', (D Ellwood, P Ozsv{\'a}th, A Stipsicz, Z Szab{\'o}, editors),
  volume~5, AMS and Clay Mathematics Institute (2006)  3--28

\bibitem{Lectures-HFH}
\textbf{P Ozsv{\'a}th}, \textbf{Z Szab{\'o}}, \emph{Lectures on {H}eegaard
  {F}loer homology}, from: ``Floer {H}omology, {G}auge {T}heory, and {L}ow
  {D}imensional {T}opology'', (D Ellwood, P Ozsv{\'a}th, A Stipsicz, Z
  Szab{\'o}, editors), volume~5, AMS and Clay Mathematics Institute (2006)
  29--70

\bibitem{Rasmussen-Knot}
\textbf{J Rasmussen}, \emph{Floer homology and knot complement}, PhD thesis,
  Harvard University (2003)Ph. D. thesis

\bibitem{Lens-Rasmussen}
\textbf{J Rasmussen}, \emph{Lens space surgeries and {L}-space homology
  spheres} (2007) \xox{arXiv}{0710.2531v1}

\bibitem{Rong}
\textbf{Y\,W Rong}, \emph{Some Knots Not Determined by Their Complements},
  from: ``Quantum Topology'', Ser. on Knots and Everything (Book 3), World
  Scientific Pub Co Inc (1993)  339--353

\bibitem{Wu}
\textbf{Z Wu}, \emph{Cosmetic surgery in {L}-space homology spheres}, Geom.
  Topol. 15 (2011) 1157--1168

\end{thebibliography}

\end{document}